

\baselineskip=14pt
\parskip=10pt

\font\eightrm=cmr8 
\font\eighttt=cmtt8
\magnification=\magstephalf
\def\P{{\cal P}}

\def\1{{\overline{1}}}
\def\2{{\overline{2}}}
\parindent=0pt
\overfullrule=0in
\def\Tilde{\char126\relax}
\def\frac#1#2{{#1 \over #2}}
\bf
\centerline
{
The Discrete Analog of the Malgrange-Ehrenpreis Theorem
}
\rm
\bigskip
\centerline{ {\it
Doron 
ZEILBERGER}\footnote{$^1$}
{\eightrm  \raggedright
Department of Mathematics, Rutgers University (New Brunswick),
Hill Center-Busch Campus, 110 Frelinghuysen Rd., Piscataway,
NJ 08854-8019, USA.
{\eighttt zeilberg  at math dot rutgers dot edu} ,
\hfill \break
{\eighttt http://www.math.rutgers.edu/\~{}zeilberg/} .
First version: July 21, 2011.
Accompanied by Maple package {\eighttt LEON} available from
{\eighttt http://www.math.rutgers.edu/\~{}zeilberg/tokhniot/LEON} .
Supported in part by the USA National Science Foundation.
}
}

\qquad \qquad \qquad \qquad \qquad \qquad  \qquad \qquad \qquad \qquad \qquad \qquad 
{\it In fond meomory of Leon Ehrenpreis}

One of the {\bf landmarks} of the modern theory of partial {\bf differential} equations is the
{\bf Malgrange-Ehrenpreis}[E][M] theorem (see [Wi]) that states that
every non-zero linear partial differential operator with constant coefficients has a Green's function
(alias {\bf fundamental solution}). Recently Wagner[W] gave an elegant {\it constructive} proof.

In this short note I will state the {\bf discrete} analog,
and give two proofs. The first
one is Ehrenpreis-style, using {\it duality}, and the second one is {\it constructive},
using {\it formal Laurent series}. 

Let $Z$ be the set of integers, and $n$ a positive integer. Consider functions $f(m_1, \dots, m_n)$
from $Z^n$ to the complex numbers (or any field).
A {\bf linear partial difference operator with constant coefficients} $\P$ is anything of the form
$$
\P f(m_1, \dots, m_n):=
\sum_{\alpha \in A} c_{\alpha} f(m_1+\alpha_1, \dots, m_n+\alpha_n) \quad,
$$
where $A$ is a {\bf finite} subset of $Z^n$ and $\alpha=(\alpha_1, \dots, \alpha_n)$,
and the $c_\alpha$  are {\bf constants}.

For example, the {\bf discrete Laplace operator} in two dimensions:
$$
f(m_1,m_2) \rightarrow f(m_1,m_2)-\frac{1}{4} ( f(m_1+1,m_2)+f(m_1-1,m_2)+f(m_1,m_2+1)+f(m_1,m_2-1) ) \quad .
$$

The {\bf symbol} of the operator $\P$ is the {\bf Laurent polynomial}
$$
P(z_1, \dots, z_n)=\sum_{\alpha \in A} c_{\alpha} z_1^{\alpha_1} \cdots, z_n^{\alpha_n} \quad.
$$

The {\bf discrete delta function} is defined in the obvious way
$$
\delta(m_1, \dots, m_n)=
\cases{ 1,&if $(m_1, \dots, m_n) =(0,0, \dots , 0)$;\cr
0,& otherwise.}
$$
Note that the beauty of the discrete world is that the delta function is a {\it genuine} function, not
a ``generalized'' one, and one does not need the intimidating edifice of Schwartzian {\it distributions}

We are now ready to state the

{\bf Discrete Malgrange-Ehrenpreis Theorem}: Let $\P$ be any non-zero {\bf linear partial difference operator
with constant coefficients}. There exists a function $f(m_1, \dots, m_n)$ defined on $Z^n$ such that
$$
\P f(m_1, \dots, m_n)=\delta(m_1, \dots , m_n) \quad .
$$

{\bf First Proof (Ehrenpreis-style)} Consider the infinite-dimensional vector space,
$C[z_1, \dots, z_n,z_1^{-1}, \dots, z_n^{-1}]$,
of {\it all} Laurent polynomials
in $z_1, \dots, z_n$. Every function $f$ on $Z^n$ uniquely defines a {\bf linear functional} $T_f$ defined
on monomials by
$$
T_f(z_1^{m_1} \cdots z_n^{m_n}):=f(m_1, \dots, m_n) \quad,
$$
and extended by linearity. Conversely, any linear functional gives rise to a function on $Z^n$.
Let $P(z_1, \dots, z_n)$ be the symbol of the operator $\P$.
We are looking for a linear functional $T$ such that for every monomial $z_1^{m_1} \cdots z_n^{m_n}$
$$
T(P(z_1, \dots z_m) z_1^{m_1} \cdots z_n^{m_n})=T_{\delta}( z_1^{m_1} \cdots z_n^{m_n} ) \quad .
$$
By linearlity,  for {\it any} Laurent polynomial $a(z_1, \dots, z_n)$
$$
T(P(z_1, \dots, z_n) a(z_1, \dots, z_n)=T_{\delta}( a(z_1, \dots, z_n) ) \quad .
$$
So $T$ is defined on the (vector) {\bf subspace} $P(z_1, \dots, z_n)C[z_1, \dots, z_n,z_1^{-1}, \dots, z_n^{-1}]$
of $C[z_1, \dots, z_n,z_1^{-1}, \dots, z_n^{-1}]$.
By elementary linear algebra, every linear functional on  the former 
can be {\bf extended} (in many ways!) to the latter. {\bf QED}.

Before embarking on the second proof we have to recall the notion of {\bf formal power series} and 
more generally {\bf formal Laurent series}.

A {\it formal power series} in one variable $z$ is any creature of the form
$$
\sum_{i=0}^{\infty} a_i z^i \quad.
$$
More generally, 
a {\bf positive} {\it formal Laurent series} is any creature of the form
$$
\sum_{i \geq m}^{\infty} a_i z^i \quad,
$$
where $m$ is a (possibly negative) integer.
On the other hand
a {\bf negative} {\it  formal Laurent series} is any creature of the form
$$
\sum_{i \leq m}^{\infty} a_i z^i \quad,
$$
where $m$ is a (possibly positive) integer.

A {\bf bilateral} formal Laurent series goes {\bf both ways}
$$
\sum_{i=-\infty}^{\infty} a_i z^i \quad.
$$

Note that the class of bilateral formal Laurent series is an abelian additive group, but one {\bf can't}  multiply there.
On the other hand one can legally multiply two positive formal Laurent series by each other, and 
two negative formal Laurent series by each other, but don't mix them!
Of course it is always legal to multiply any Laurent polynomial by any  bilateral formal power series.
{\bf But watch out for} {\it zero-divisors}, e.g.
$$
(1-z) \sum_{i=-\infty}^{\infty} z^i \, = \, 0 \quad .
$$

Any Laurent polynomial $p(z)=a_iz^i+ \dots a_jz^j$ of low-degree $i$ and (high-)degree $j$ in  $z$ 
(so $a_i \neq 0$, $a_j \neq 0$)
has
{\bf two} multiplicative inverses. One in the ring of positive Laurent polynomials, and the other
in the ring of negative Laurent polynomials. Simply  write $p(z)=z^i a_i p_0(z)$ and
get $1/p(z)=z^{-i} (1/a_i) p_0(z)^{-1}$, and writing $p_0(z)=1+ q_0(z)$, we form
$$
p_0(z)^{-1}=(1+q_0(z))^{-1}=\sum_{i=0}^{\infty} (-1)^i q_0(z)^i \quad ,
$$
and this makes perfect sense and {\it converges} in the ring of formal power series.
Analogously one can form a multiplicative inverse in powers in $z^{-1}$.

It follows that every {\it rational function} $P(z)/Q(z)$ in one variable, $z$,  has (at least) two inverses, one pointing
positively, one negatively.

What about a rational function of several variables, $P(z_1, \dots, z_n)/Q(z_1, \dots ,z_n)$?
Here we can form $2^nn!$ inverses. There are $n!$ ways to order the variables, and for each
of these one can decide whether to do the postive-pointing inverse or the negative-pointing one.
At each stage we get a formal one-sided formal Laurent series whose coefficients are rational
functions of the remaining variables, and one just keeps going.

{\bf Second Proof (Constructive)}: To every discrete function $f(m_1, \dots, m_n)$ associate the
bilateral formal Laurent series
$$
\sum_{(m_1, \dots, m_n) \in Z^n} f(m_1, \dots, m_n) z_1^{m_1} \cdots z_n^{m_n} \quad.
$$
We need to ``solve'' The equation
$$
P(z_1^{-1}, \dots, z_n^{-1}) \left ( \sum_{(m_1, \dots, m_n) \in Z^n} f(m_1, \dots, m_n) z_1^{m_1} \cdots z_n^{m_n} \right )=1 \quad .
$$
So ``explicitly''
$$
 \sum_{(m_1, \dots, m_n) \in Z^n} f(m_1, \dots, m_n) z_1^{m_1} \cdots z_n^{m_n} =1/P(z_1^{-1}, \dots , z_n^{-1}) \quad ,
$$
and we just described how to do it in $2^nn!$ ways.

{\bf The Maple package LEON}

This article is accompanied by a Maple  package {\tt LEON}.
One of its numerous procedures is {\tt FS}, that implements the above constructive proof. {\tt LEON} can also compute
{\bf polynomial bases} to solutions of linear partial difference equations with constant coefficients, 
compute {\it Hilbert Series} for spaces of solutions of systems of linear differential equations,
as well as
find {\it multiplicity varieties}, in the style of Ehrenpreis, to $0$-dimensional ones.

{\bf Leon Ehrenpreis (1930-2010) A truly FUNDAMENTAL Mathematician (a Videotaped lecture)}

I strongly urge readers to watch my lecture, available in six parts from {\tt YouTube}, and
in two parts from {\tt Vimeo}, see:

{\tt http://www.math.rutgers.edu/\Tilde zeilberg/mamarim/mamarimhtml/leon.html} \quad  .

That page contains links to both versions, as well as numerous input and output files
for the Maple package {\tt LEON}.

{\bf References}

[E] Leon Ehrenpreis, {\it Solution of some problems of division. I. Division by a polynomial of derivation}, Amer. J. Math. {\bf 76}(1954), 883-903.

[M] Bernard Malgrange, 
{\it Existence et approximation des solutions des \'equations aux d\'eriv\'es partielles et des \'equations de convolution}, 
Ann. Inst. Fourier, Grenoble {\bf 6}(1955-1956): 271-355.

[W] Peter Wagner, {\it A new constructive proof of the Malgrange-Ehrenpreis theorem}, Amer. Math. Monthly {\bf 116}(2009), 457-462.

[Wi] Wikipedia,the free Encyclopedia, {\it Malgrange-Ehrenpreis Theorem}, 
Retrieved 16:10, July 21, 2011,

\end